\documentclass[11pt,a4paper]{amsart}

\usepackage[T1]{fontenc}
\usepackage{lmodern}
\usepackage[margin=24mm]{geometry}
\usepackage{amsmath,amssymb,amsthm,mathtools}
\usepackage{enumitem}
\usepackage{tikz}
\usepackage{microtype}
\usepackage{comment}

\numberwithin{equation}{section}
\newtheoremstyle{sourceplain}{8pt}{8pt}{\normalfont}{}{\bfseries}{}{0.6em}{}
\theoremstyle{sourceplain}
\newtheorem{theorem}{Theorem}[section]
\newtheorem{proposition}[theorem]{Proposition}
\newtheorem{lemma}[theorem]{Lemma}

\setlist[enumerate]{topsep=4pt,itemsep=4pt,parsep=2pt}
\title{Local Square-Root Cancellation for Quartic Weyl Sums}

\author{Guanyu Fang}
\address{Guanyu Fang, School of Mathematics and Statistics, Xi'an Jiaotong University, Xi'an, Shaanxi, China.}

\email{guanyu.fang@foxmail.com}

\author{Sian Fang}
\address{Sian Fang, School of Mathematical Sciences,
Fudan University, People's Republic of China}
\email{24110180013@m.fudan.edu.cn}

\author{Victor Zhenyu Guo}
\address{Victor Zhenyu Guo, School of Mathematics and Statistics, Xi'an Jiaotong University, Xi'an, Shaanxi, China.}

\email{vzguo@foxmail.com; guozyv@xjtu.edu.cn}

\author{Danqing He}
\address{Danqing He, School of Mathematical Sciences,
Fudan University, People's Republic of China}
\email{hedanqing@fudan.edu.cn}

\thanks{D. He is supported by  National Key R$\&$D Program of China (No. 2021YFA1002500),  NNSF of China (No. 12322105), and the New Cornerstone Science Foundation. V. Z. Guo is supported by the National Natural Science Foundation of China (No. 12271422)}

\date{}
\begin{document}
\maketitle

\begin{abstract}
We study local $L^p$ moments of quartic Weyl sums along $(t,t^2,t^3,t^4)$ over subintervals of $[N/2,N]$ on positive boxes $[0,N]\times[0,N^2]\times[0,N^a]\times[0,N^b]$. For $11\le p\le 12$, $a\ge b\ge 0$, and $a+b=p/2-3$, we prove $L^p$-averaged square-root cancellation up to an $N^\varepsilon$ loss, establishing the quartic moment-curve case of Demeter's conjecture on positive boxes. For $p=10$, $1\le a\le 2$, and $b=2-a$, we establish the corresponding tenth-moment bound uniformly in $a$, even though square-root cancellation fails on symmetric boxes when $b>0$. The proofs combine Poisson summation on rational arcs, curvature decomposition and rescaling, and small-cap decoupling for the quartic and cubic moment curves.
\end{abstract}


\section{Introduction}

Let $\mathrm{e}(t)=\exp(2\pi it)$. For integers $k\ge2$ and $N\ge2$, the
classical polynomial Weyl sum is
\begin{equation}
 W_k(\boldsymbol\alpha;N)
 =\sum_{n=1}^N\mathrm{e}(\alpha_1n+\cdots+\alpha_kn^k),
 \qquad \boldsymbol\alpha\in\mathbb R^k.
 \label{eq:classical-weyl}
\end{equation}
The study of this subject grew
from Weyl's work on uniform distribution modulo one.
The Weyl sum conjecture is as follows: if $(r,q)=1$,
$|\alpha_k-r/q|<q^{-2}$, and $N\le q\le N^{k-1}$, then for any $\epsilon>0$, there exists a constant $C_{k,\epsilon}>0$ such that
\begin{equation}
 \sup_{\alpha_1,\ldots,\alpha_{k-1}\in\mathbb R}
 |W_k(\boldsymbol\alpha;N)|
 \le C_{k,\epsilon} N^{1+\epsilon}
 \left(q^{-1}+\frac{q}{N^k}\right)^{1/k}.
 \label{eq:pointwise-conjecture}
\end{equation}
The classical Weyl differencing yields the estimate
\begin{equation}
 |W_k(\boldsymbol\alpha;N)|
 \ll_{k,\epsilon }N^{1+\epsilon }
 \left(q^{-1}+N^{-1}+qN^{-k}\right)^{2^{1-k}}.
 \label{eq:classical-weyl-bound}
\end{equation}
For $k\ge3$, the modern Vinogradov mean-value theorem \cite{BDG16}  gives the exponent   $1/[k(k-1)]$.

In this note, we study a special Weyl sum which is related to the Lindel\"of hypothesis \cite{Bourgain}.
More precisely,
let $N\ge2$ be an integer and put
\[
\Gamma(t)=(t,t^2,t^3,t^4).
\]
For an interval $I\subset[N/2,N]$, we define
\begin{align}
	S_{I,N}(x)&=\sum_{n\in I\cap\mathbb Z}\mathrm{e}(x\cdot\Gamma(n/N)).
	\label{eq:sum-definition}
\end{align}
Instead of considering its pointwise estimate as in \eqref{eq:pointwise-conjecture}, we study a local mean estimate.
For real parameters $a,b\ge0$, denote
\begin{align}
	\Omega_N(a,b)&=[0,N]\times[0,N^2]\times[0,N^a]\times[0,N^b].
	\label{eq:box-definition}
\end{align}
Our main result is as follows.

\begin{theorem}\label{thm:main}

	Let $11\le p\le12$, $a\ge b\ge0$, and $a+b=\frac p2-3$.
Then, for any $\varepsilon>0$, there exists $ C_\varepsilon>0$ such that
	\begin{equation}\label{eq:main-moment}
		\int_{\Omega_N(a,b)}|S_{I,N}(x)|^p\,dx
		\lesssim_{p,a,\varepsilon}N^{p+\varepsilon}.
	\end{equation}
for all $ N\ge2$ and $I\subset[\frac N2,N]$.

\end{theorem}

The Lindel\"of hypothesis concerns the growth of the zeta function on the critical line, which asserts that
$|\zeta(1/2+it)|\ll_\varepsilon (1+|t|)^\varepsilon$.
Bourgain \cite{Bourgain} proved that $|\zeta(1/2+it)|\lesssim_\varepsilon (1+|t|)^{13/84+\varepsilon}$, where a related curve $
(t,t^2,t^{3/2},t^{1/2})$, proposed by Huxley \cite{Huxley96}, was used to prove the first spacing problem in the Bombieri–Iwaniec method. Bourgain's key estimate actually gives also \eqref{eq:main-moment} for $(p,a,b)=(12,2,1)$, which was later stated explicitly by Demeter \cite{Demeter}. Demeter~\cite{Demeter} proved also the case $p=12$, $3/2\le a\le2$, $b=3-a$. Glidewell \cite{Glidewell} obtained the range
$11\le p\le14$, $2\le a\le3$, $0<b\le1$,
$a+b=p/2-3$.

Since $|\Omega_N(a,b)|=N^{3+a+b}=N^{p/2}$, Theorem~\ref{thm:main} is
equivalent to
\[
\left(\frac1{|\Omega_N(a,b)|}
\int_{\Omega_N(a,b)}|S_{I,N}(x)|^p\,dx\right)^{1/p}
\lesssim_{p,a,\varepsilon}N^{1/2+\varepsilon},
\]
expressing square-root cancellation on average over $\Omega_N(a,b)$
(after rescaling $\varepsilon$ by $p$). Under $x_1=Ny_1$, $x_2=N^2y_2$
(Jacobian $N^3$), Theorem~\ref{thm:main} establishes the quartic
moment-curve case of \cite[Conjecture 8.1]{Demeter} on the positive
boxes \eqref{eq:box-definition}.
For \(p<11\) and \(b>0\), where Demeter showed that square-root
cancellation fails on symmetric boxes \cite[Theorem~8.2]{Demeter},
the positive-box bounds for the unweighted
sum appear new both for \(10<p<11,\ (p-6)/4\le a\le 2(p-6)/5\)
(Proposition~\ref{prop:QuarticWeylSum:common-range}) and for
\(p=10,\ 1\le a<2\) (Theorem~\ref{thm:tenth}).

Concerning the case $p=10$, we have the following result.

\begin{theorem}\label{thm:tenth}

Let $1\le a\le2$ and $b=2-a$. Then for every $\varepsilon>0$ there exists a constant $C_\varepsilon>0$ such that
\begin{equation}\label{eq:tenth-upper}
 \int_{\Omega_N(a,b)}|S_{I,N}(x)|^{10}\,dx
       \lesssim_{\varepsilon} N^{10+\varepsilon}
\end{equation}
for all 
$a$, $N\ge2$, and $I\subset[N/2,N]$,

\end{theorem}

Bourgain's earlier tenth-moment bound for the same monomials loses a
power of \(N\) when \(b>0\) \cite[Theorem~4, (4.18) and
(5.20)--(5.21)]{BourgainMeanValues}, and the endpoint
\((p,a,b)=(10,2,0)\) follows by fixing \(x_4\) and applying
\cite[Theorem~3.3]{DemeterGuthWang}.

\medskip

The paper is organized as follows. Section 2 passes to a smoothed sum and combines Poisson summation with oscillatory-integral estimates to split off and bound the near terms. Section 3 decomposes the remainder into dyadic curvature bands, establishes short-interval pointwise bounds, and rescales each band to reduce its $L^p$ integral to localized moments of a rescaled sum. Section 4 uses quartic small-cap decoupling to derive localized $L^p$ estimates $(10<p\le12)$, proves Proposition 4.2, and combines these estimates with earlier results to establish Theorem 1.1.
Section 5 proves Theorem 1.2 by fixing the fourth rescaled coordinate and applying cubic small-cap decoupling.

\medskip

\noindent{\bf Notations.}
We denote $
e(x):=\exp(2\pi ix).$
By
$
A\lesssim B
$
we mean that there exists an absolute constant $C>0$ such that $A\le CB$.
If additionally $C$ depends on $\epsilon$, we also write
$
A\lesssim_\epsilon B.$
We denote $A\sim B$ when $A\lesssim B$ and $B\lesssim A$.
By $A\ll B$ we mean that $A$ is much smaller than $B$, depending on how small we need in a specific situation.
When $A\lesssim_{\varepsilon}T^{\varepsilon}B$ for all $\varepsilon>0$ we denote
$
A\lessapprox B.
$

\medskip

\noindent{\bf Use of AI.} The authors identified the research direction and led the work, while AI tools supplied some essential arguments and ideas. The authors independently verified every proof detail, rewrote the manuscript, and take full responsibility for its content.

\section{Poisson summation and near terms}
\label{sec:poisson}

Fix a smooth cutoff $0\le\gamma\le1$ in $C_c^\infty((1/4,5/4))$ equal
to $1$ on $[1/2,1]$, and define
\begin{align*}
	I(A,B,u,v)&=\int_{\mathbb R}\gamma(t)\mathrm{e}(At+Bt^2+ut^3+vt^4)\,dt,\\
	S_\gamma(x)&=\sum_{n\in\mathbb Z}\gamma(n/N)\mathrm{e}(x\cdot\Gamma(n/N)).
\end{align*}

We give the following initial lemma. 

\begin{lemma}\label{lem:continuous}
	For all $u,v\ge0$, put $T=1+u+v$.
	\begin{enumerate}
		\item Uniformly in $u,v\ge0$,
		\begin{equation}\label{eq:continuous-tenth}
			\int_{\mathbb R^2}|I(A,B,u,v)|^{10}\,dA\,dB\lesssim T^{-2},
		\end{equation}
		and hence for all $U,V\ge0$,
		\begin{equation}\label{eq:tenth-continuous-box}
			\int_{\mathbb R^2\times[0,U]\times[0,V]}|I(A,B,u,v)|^{10}\,dA\,dB\,du\,dv
			\lesssim\log\bigl(1+\min(U,V)\bigr).
		\end{equation}
		\item For every $p>10$, with $\mathbb R_+=[0,\infty)$,
		\[
		\int_{\mathbb R^2\times\mathbb R_+^2}|I(A,B,u,v)|^p\,dA\,dB\,du\,dv
		\lesssim_p 1.
		\]
	\end{enumerate}
\end{lemma}
\begin{proof}
	
We write $P(t) = Bt^2 + ut^3 + vt^4$ for short and take $C = 25$. Firstly, we give the proof of equation \eqref{eq:continuous-tenth}. Considering the case when $|B|>CT$ and any $T\ge1$, it is clear that
$$
|P''(t)|=|2B+6ut+12vt^2|\asymp|B|
$$
on $[1/4,5/4]$. A second-derivative test (see \cite[Lemmas 2.5]{TaoStationary}) gives that $\|I(\cdot,B,u,v)\|_\infty\lesssim|B|^{-1/2}$. We apply the Plancherel Theorem on $A$ to get that 
$$
\|I(\cdot,B,u,v)\|_{L^2(\mathbb R)}^2=\|\gamma\|_2^2\lesssim1.
$$
Hence
\[
\int_{\mathbb R}|I(A,B,u,v)|^{10}\,dA
\le\|I(\cdot,B,u,v)\|_\infty^8\|I(\cdot,B,u,v)\|_2^2
\lesssim|B|^{-4},
\]
and integrating over $|B|>CT$ yields $O(T^{-3})\le O(T^{-2})$.

Next we work on the case when $|B|\le CT$. If $1\le T<2$ we have $\|I\|_\infty\le1$ and $\|I\|_2^2\lesssim1$ which give that $\int_{\mathbb R}|I|^{10}\,dA\lesssim1$. Integrating over $|B|\le2C$ yields the bound $O(1)=O(T^{-2})$. Now for $T \geq 2$, it follows that $u+v=T-1\ge T/2$ and $P'''(t)=6u+24vt\asymp T$ on $[1/4,5/4]$. We partion the interval $[1/4,5/4]$ into the interval
$$
E_0=\{t\in[1/4,5/4]:|P''(t)|\le T^{2/3}\}
$$
and dyadic intervals 
$$
E_{\mu,\pm}=\{t\in[1/4,5/4]:\mu<\pm P''(t)\le2\mu\}
$$
for $\mu=2^rT^{2/3}\lesssim T$. Note that the length of $E_0$ and $E_{\mu,\pm}$ are $O(T^{-1/3})$ and $O(\mu/T)$ respectively. We define that
$$
I_E(A)=\int_E\gamma(t)\mathrm{e}(At+P(t))\,dt.
$$
With Plancherel Theorem we have that $\|I_E\|_2^2\le|E|$. Moreover also with $\|I_{E_0}\|_\infty\le|E_0|$ it gives that $\|I_{E_0}\|_{10}^{10}\le|E_0|^9\lesssim T^{-3}$. On each $E_{\mu,\pm}$,
the second-derivative test gives $\|I_{E_{\mu,\pm}}\|_\infty\lesssim\mu^{-1/2}$,
so
\[
\|I_{E_{\mu,\pm}}\|_{10}^{10}
\le\|I_{E_{\mu,\pm}}\|_\infty^8\|I_{E_{\mu,\pm}}\|_2^2
\lesssim\mu^{-4}\cdot\frac{\mu}{T}=T^{-1}\mu^{-3}.
\]
Summing the tenth roots by Minkowski's inequality gives
$\int_{\mathbb R}|I|^{10}\,dA\lesssim T^{-3}$ uniformly for $|B|\le CT$.
Integrating over $|B|\le CT$ gives $O(T^{-2})$, proving
\eqref{eq:continuous-tenth}. Integrating \eqref{eq:continuous-tenth}
over $[0,U]\times[0,V]$ yields
\[
\int_0^U\int_0^V\frac{du\,dv}{(1+u+v)^2}
=\log\frac{(1+U)(1+V)}{1+U+V}
\le\log\bigl(1+\min(U,V)\bigr),
\]
proving \eqref{eq:tenth-continuous-box}.

Finally we consider the case when $p > 10$. On $[1/4,5/4]$,  we have
$$
P'''(t)=6u+24vt\asymp u+v,
$$
so the third-derivative test (\cite[Lemmas 2.5]{TaoStationary}) together with $|I|\le1$ yields $\|I(\cdot,\cdot,u,v)\|_\infty\lesssim T^{-1/3}$. Combining this with
\eqref{eq:continuous-tenth} gives
\[
\int_{\mathbb R^2}|I(A,B,u,v)|^p\,dA\,dB
\le\|I(\cdot,\cdot,u,v)\|_\infty^{p-10}
\int_{\mathbb R^2}|I(A,B,u,v)|^{10}\,dA\,dB
\lesssim T^{-2-(p-10)/3}.
\]
Since $2+(p-10)/3>2$, integrating over $(u,v)\in\mathbb R_+^2$ gives $O(1)$.
\end{proof}	

Now we write the phase point in the form $x=(A,N^2c/q+B,u,v)$, where $q\ge1$, $(c,q)=1$. Throughout the paper, put
\[
P(t)=Bt^2+ut^3+vt^4,\qquad \lambda=\frac Nq,\qquad
G_q(m,c)=\sum_{r\bmod q}\mathrm{e}\Bigl(\frac{cr^2+mr}{q}\Bigr).
\]
For any $w\in C_c^\infty(\mathbb R)$, we apply Poisson summation separately to each residue class $n\equiv r\bmod q$. Write $n =r + qk$. 
\begin{align}\label{eq:poisson}
	\nonumber &\qquad \sum_{n\in\mathbb Z} w\big( \frac{n}{N} \big) \mathrm{e}\bigl(\frac{cn^2}{q}+\frac{An}{N}+P(\frac{n}{N})\bigr) \\ 
	\nonumber &= \sum_{r \bmod q}  \mathrm{e} \big(\frac{cr^2}{q} \big) \sum_{k \in \mathbb{Z}} w( \frac{r + qk}{N} ) \mathrm{e} \big(A \frac{r+qk}{N} + P (\frac{r+qk}{N} ) \big) \\
	&= \frac{N}{q} \sum_{m \in \mathbb{Z}} G_q(m,c) I_m(w),
\end{align}	
where
$$
I_m(w)=\int_{\mathbb R}w(t)\mathrm{e}((A-m\lambda)t+P(t))\,dt.
$$
In particular, taking $w=\gamma$ gives
$$
S_\gamma(A,N^2c/q+B,u,v)=\frac Nq\sum_m G_q(m,c)I(A-m\lambda,B,u,v).
$$
For all $q\ge1$, the quadratic Gauss sum satisfies
\[
|G_q(m,c)|^2
=\sum_{h\bmod q}\mathrm{e}\Bigl(\frac{ch^2+mh}{q}\Bigr)
\sum_{r\bmod q}\mathrm{e}\Bigl(\frac{2chr}{q}\Bigr)
\le q\gcd(2,q)\le 2q,
\]
so $|G_q(m,c)|\le\sqrt{2q}$.

Huxley and Kolesnik have already used Poisson summation with a rational quadratic
phase and a partition by curvature in their work on the first spacing
problem \cite{HuxleyKolesnikIII}.

\begin{lemma}\label{lem:nonstationary}
Let $B,u,v\in\mathbb R$, $1\le q\le N$, $(c,q)=1$, and $\lambda=N/q$. Let $I_w\subset\mathbb R$ be an interval of length $O(\Delta)$ with $\Delta>0$. Let $w\in C_c^\infty(I_w)$ and $B_w=\max_{0\le r\le 2}\Delta^r\|w^{(r)}\|_\infty$. Assume that $\Delta|P''|\lesssim\lambda$ and $\Delta^2|P'''|\lesssim\lambda$ on $I_w$. Then, uniformly in $A\in\mathbb R$,
\begin{equation}\label{eq:nonstationary}
\frac Nq\sum_{\substack{m\in\mathbb Z\\ \operatorname{dist}(m\lambda-A,\,P'(I_w))\ge\lambda}} |G_q(m,c)|\,|I_m(w)| \lesssim \frac{\sqrt q\,B_w}{\lambda\Delta},
\end{equation}
where $dist(t, \mathcal{T})$ stands the distance between a point $t$ and a set $\mathcal{T}$.  
\end{lemma}

\begin{proof}
Fix $t_*\in I_w$. Write
$$
\tilde w(s) = w(t) = w(t_*+\Delta s)
$$
and
$$
\Phi_m(s)=(A-m\lambda)(t_*+\Delta s)+P(t_*+\Delta s). 
$$
Then 
$$
I_m(w)=\Delta\int\tilde w(s)\mathrm{e}(\Phi_m(s))\,ds.
$$
Also easy to check that $\operatorname{supp}\tilde w$ has length $O(1)$, 
$\|\tilde w^{(r)}\|_\infty\le B_w$ ($0\le r\le2$), and
$|\Phi_m''(s)|,|\Phi_m'''(s)|\lesssim Z$ on $\operatorname{supp}\tilde w$, where $Z = \lambda \Delta$. 

If $\operatorname{dist}(m\lambda-A,P'(I_w))\in[j\lambda,(j+1)\lambda)$
for an integer $j\ge1$, we have that $|\Phi_m'(s)|\ge jZ$ on
$\operatorname{supp}\tilde w$. Integrating by parts twice gives that
\[
\int_{\mathbb R}\tilde w(s)\mathrm{e}(\Phi_m(s))\,ds
=-\frac{1}{4\pi^2}\int_{\mathbb R}
\left(\frac{\tilde w''}{(\Phi_m')^2}
-\frac{3\tilde w'\Phi_m''+\tilde w\Phi_m'''}{(\Phi_m')^3}
+\frac{3\tilde w(\Phi_m'')^2}{(\Phi_m')^4}\right)
\mathrm{e}(\Phi_m(s))\,ds,
\]
so 
$$
|I_m(w)|\lesssim B_w\Delta j^{-2}Z^{-2}=B_w j^{-2}\lambda^{-2}\Delta^{-1}.
$$
Since $P'(I_w)$ is an interval, at most two points of
$\lambda\mathbb Z-A$ lie at distance in $[j\lambda,(j+1)\lambda)$ from
$P'(I_w)$. Using $|G_q(m,c)|\le\sqrt{2q}$, summing
over $j\ge1$ proves \eqref{eq:nonstationary}.
\end{proof}

Fix $0<t_-<1/4<5/4<t_+$. For $u,v\ge0$, set
\[
\mathcal C=\{P'(t_-),P'(t_+)\}
\cup\{P'(t):t\in[t_-,t_+],\ P''(t)=0\}.
\]
It is clear that $|\mathcal C|\le3$. For a fixed large constant $C_E$ which is chosen later, we define an exceptional set $E$ by 
$$
 E=\{m\in\mathbb Z:\operatorname{dist}(m\lambda-A,\mathcal C)\le C_E\lambda\}. 
$$
We also split $S_\gamma$ by setting
$$
 S_{\mathrm{near}}(x)=\frac Nq\sum_{m\in E}G_q(m,c)I(A-m\lambda,B,u,v). 
$$
and
$$
\mathcal R=S_\gamma-S_{\mathrm{near}}.
$$
We have that $|E| \lesssim 1$ uniformly in $A, B, u, v, q$. For a function $f \ge 0$, we have the periodicity identity 
\begin{equation}\label{eq:pi}
\int_0^N\sum_{m\in\mathbb Z}f(A-m\lambda)\,dA=q\int_{\mathbb R}f(A)\,dA. 
\end{equation}
Applying H\"older's inequality over $m\in E$ and equation \eqref{eq:poisson}, it follows that
\begin{align*}
 |S_{\mathrm{near}}|^p &\le \lambda^p |E|^{p-1} \sum_{m \in E} |G_q(m,c)|^p |I(A-m\lambda, B, u, v)|^p \\
 &\lesssim N^p q^{-p/2} \sum_{m \in E} |I(A-m\lambda, B, u, v)|^p. 
\end{align*}
Let 
$$
\operatorname{arc}(c/q)=[0,N]\times J_{q,c}\times[0,U]\times[0,V]
$$
where $J_{q,c} \subset\mathbb R$, $U,V\ge0$. With \eqref{eq:pi} we have that
\begin{align*}
	 \int_{\operatorname{arc}(c/q)}|S_{\mathrm{near}}|^p &\lesssim N^p q^{-p/2} \int \, dBdudv \sum_{m \in \mathbb{Z}} \int_0^N |I(A-m\lambda, B, u, v)|^p \, dA \\
	 &\le N^p q^{-p/2+1} \int_{\mathbb{R} \times \mathbb{R} \times [0,U] \times [0,V]} |I(A, B, u, v)|^p \,dA dB du dv.
\end{align*}
Summing over the $\varphi(q)\le q$ reduced residues $c\bmod q$ and all $q\ge1$ gives
\[
\sum_{q\ge1}\sum_{\substack{c\bmod q\\(c,q)=1}}
\int_{\operatorname{arc}(c/q)}|S_{\mathrm{near}}|^p
\lesssim_p N^p\sum_{q\ge1}q^{2-p/2}
\int_{\mathbb R^2\times[0,U]\times[0,V]}|I|^p.
\]
Since $2-p/2\le-3$ for $p\ge10$, Lemma~\ref{lem:continuous} yields both
\begin{align}
	\sum_{q\ge1}\sum_{\substack{c\bmod q\\(c,q)=1}}
	\int_{\operatorname{arc}(c/q)}|S_{\mathrm{near}}|^p
	&\lesssim_p N^p && (p>10),\label{eq:exceptional}\\
	\sum_{q\ge1}\sum_{\substack{c\bmod q\\(c,q)=1}}
	\int_{\operatorname{arc}(c/q)}|S_{\mathrm{near}}|^{10}
	&\lesssim N^{10}\log\bigl(1+\min(U,V)\bigr) && (p=10).\label{eq:tenth-near}
\end{align}
The left part is to estimate $\mathcal{R}$, which will be done in next two sections. 

To pass from $S_\gamma$ back to $S_{I,N}$ for $I\subset[N/2,N]$ without
logarithmic loss, fix $(x_2,u,v)$ and put $\theta=A/N\in\mathbb T$.
Since $\gamma=1$ on $[1/2,1]$, sharp Fourier projection onto the
consecutive integers in $I$ maps $S_\gamma$ to $S_{I,N}$; as a
difference of two modulated Riesz projections, it is bounded on
$L^p(\mathbb T)$ uniformly in $I$ for $1<p<\infty$ by
\cite{HollenbeckVerbitsky}. Integrating over $(x_2,u,v)$ yields
\begin{equation}\label{eq:interval-projection}
	\int_{\Omega_N(a,b)}|S_{I,N}|^p\,dx
	\le C_p\int_{\Omega_N(a,b)}|S_\gamma|^p\,dx.
\end{equation}

\section{Curvature bands and rescaling} \label{sec:3}

In this section, we attempt to work on $\mathcal{R}$. The following lemma gives a decomposition of $\mathcal{R}$. 

\begin{lemma}\label{lem:decomposition}
	Fix a sufficiently large absolute constant $H_0$, and choose $C_E$
	sufficiently large in terms of $H_0$. Let $D\asymp1+|B|+u+v$ with $D\lesssim\lambda^2$,
	and $H=D/\lambda$.
	\begin{enumerate}
		\item If $H\le H_0$, then $\mathcal R=O(\sqrt N)$ uniformly in $A$.
		\item If $H>H_0$, put $\tau_0=c_{*} H^{-1/2}$ with a small fixed $c_{*}>0$.
		There exist smooth cutoffs $\psi_\tau$ depending on $B,u,v$ but not
		on $A$, indexed by dyadic $\tau=2^j\tau_0\lesssim1$ ($j\ge0$), such that on each $\operatorname{supp}\psi_\tau$:
		\begin{itemize}
			\item $\operatorname{supp}\psi_\tau$ lies in an interval of length
			$O(\tau)$ inside $(1/4,5/4)$;
			\item $P''$ has fixed sign and $|P''|\asymp D\tau$;
			\item $\|\partial_t^r\psi_\tau\|_\infty\lesssim_r\tau^{-r}$ for every
			$r\ge0$.
		\end{itemize}
		Defining 
		$$
		T_\tau(A)=\sum_{n\in\mathbb Z}\psi_\tau(B,u,v; \frac{n}{N}) \mathrm{e}\big(\frac{cn^2}{q}+A\frac{n}{N}+P(\frac{n}{N})\big), 
		$$
		we have uniformly in $A$,
		\begin{equation}\label{eq:decomposition}
			\mathcal R=\sum_\tau T_\tau+O(\sqrt N).
		\end{equation}
	\end{enumerate}
\end{lemma}
\begin{proof}
\emph{Case 1: $H\le H_0$.} On $[1/4,5/4]$, we have $|P''|,|P'''|\lesssim D\le H_0\lambda$, giving that
$$
|P'(t)-P'(t_-)|\lesssim H_0\lambda.
$$
Taking $C_E$ large gives
$$
\operatorname{dist}(m\lambda-A,P'([1/4,5/4]))\ge\lambda
$$
for all $m\notin E$. Applying Lemma~\ref{lem:nonstationary} with $w=\gamma$, $\Delta=1$	yields $|\mathcal R|\lesssim\sqrt q/\lambda\le\sqrt q \lesssim \sqrt N$.
	
\emph{Case 2: $H>H_0$.} In this case $D\asymp |B|+u+v$. Fix a small constant $c_0>0$. If $u+v\le c_0D$, then $|B|\asymp D$ and $|P''|\asymp D$ with fixed sign on $[t_-,t_+]$. If $u+v>c_0D$, then $P'''(t)=6u+24vt\asymp D$ on $[t_-,t_+]$; if $P''$ has no zero there, $|P''|\asymp D$ on $\operatorname{supp}\gamma$.

Otherwise $P'''\asymp D$ on $[t_-,t_+]$ and $P''$ has a unique zero $t_c\in[t_-,t_+]$, so $|P''(t)|\asymp D|t-t_c|$ on $[t_-,t_+]$. Fix
$\varphi \in C_c^\infty((-2,2))$ equal to $1$ on $[-1,1]$. 
Also define 
$$
\psi_{\mathrm{core}} = \gamma(t) \varphi(\frac{t-t_c}{\tau_0}).  
$$
It follows that
$$
\psi_{\tau,+}(t) =  \gamma(t) \big( \varphi(\frac{t-t_c}{2\tau}) - \varphi(\frac{t-t_c}{\tau}) \big) \mathbf{1}_{t>t_c},
$$
$$
\psi_{\tau,-}(t) =  \gamma(t) \big( \varphi(\frac{t-t_c}{2\tau}) - \varphi(\frac{t-t_c}{\tau}) \big) \mathbf{1}_{t<t_c}. 
$$

Hence 
$$
\gamma=\psi_{\mathrm{core}}+\sum_\tau(\psi_{\tau,+}+\psi_{\tau,-})
$$
with $\tau_0\le\tau\lesssim1$ and $\|\partial_t^r\psi_\tau\|_\infty\lesssim_r\tau^{-r}$. We note that in the other two subcases the only choice of $\psi_\tau = \gamma$ and $\psi_{\mathrm{core}} = 0$. 

Applying \eqref{eq:poisson} to $\gamma$ and each $\psi_\tau$ yields
\[
\mathcal R-\sum_\tau T_\tau =\frac Nq\sum_{m\notin E}G_q(m,c)I_m(\psi_{\mathrm{core}}) - \sum_\tau\frac Nq\sum_{m\in E}G_q(m,c)I_m(\psi_\tau).
\]
If $\psi_{\mathrm{core}}\not\equiv0$, on
$$
I_{\mathrm{core}}=[t_c-2\tau_0,t_c+2\tau_0]\cap[1/4,5/4]
$$
we have $|P''|\lesssim D\tau_0$, $|P'''|\asymp D$, and $|P'(t)-P'(t_c)|\lesssim D\tau_0^2=c^2\lambda$. By Lemma~\ref{lem:nonstationary} with $\Delta=\tau_0$ it gives that 
$$
\frac Nq\sum_{m\notin E}G_q(m,c)I_m(\psi_{\mathrm{core}}) \lesssim \frac{q^{1/2}}{\lambda \tau_0} =O(\sqrt q)\le O(\sqrt N). 
$$
Moreoever the seond-derivative test gives 
\[
\sum_\tau\frac Nq\sum_{m\in E}|G_q(m,c)|\,|I_m(\psi_\tau)|
\lesssim\frac{N}{\sqrt q}\sum_{\tau}(D\tau)^{-1/2}
\lesssim\frac{N}{\sqrt q}(D\tau_0)^{-1/2}
\asymp\sqrt N H^{-1/4}\lesssim\sqrt N. 
\]
\end{proof}

Let $a_n=\sum_j d_j\mathrm{e}(n\theta_j)$ with $C = \sum_j|d_j|$. For $1\le p\le\infty$, Minkowski's inequality gives that
\begin{equation}\label{eq:completion}
\|\sum_n b_n a_n \mathrm{e}(\frac{An}{N})\|_{L^p(0,N)}\le C \|\sum_n b_n\mathrm{e}(\frac{An}{N})\|_{L^p(0,N)}.
\end{equation}

\begin{lemma}
	\label{lem:cap}
Let $D>0$ and $0<\tau\lesssim1$ satisfy $D\tau\lesssim\lambda^2$ and $J = D\tau^2/\lambda \gtrsim 1$. Let $\Delta = \lambda/(D\tau)$. Suppose $I_{\mathrm{reg}}\subset\mathbb R$ is an interval on which
$|P''|\asymp D\tau$ and $|P'''|\lesssim D$.
\begin{enumerate}
	\item If $w\in C_c^\infty(I_{\mathrm{reg}})$ is
	supported in a subinterval of length $O(\Delta)$ and let
	$$
	B_w=\max_{0\le r\le 2}\Delta^r\|w^{(r)}\|_\infty,
	$$
	then uniformly in $A\in\mathbb R$,
	\begin{equation}\label{eq:smooth-cap}
		\sum_{n\in\mathbb Z} w\big( \frac{n}{N} \big) \mathrm{e}\bigl(\frac{cn^2}{q}+\frac{An}{N}+P(\frac{n}{N})\bigr)
		\lesssim N (qD\tau)^{-1/2} B_w.
	\end{equation}
	\item  If $\eta\in C_c^\infty(I_{\mathrm{reg}})$
	satisfies 
	$$
	\max_{0\le r\le2}\tau^r\|\eta^{(r)}\|_\infty\lesssim1,
	$$
	then for every interval $I$ of length $O(\Delta)$, uniformly in $A\in\mathbb R$,
	\begin{equation}\label{eq:sharp-cap}
		\sum_{\frac{n}{N} \in I} \eta\big( \frac{n}{N} \big) \mathrm{e}\bigl(\frac{cn^2}{q}+\frac{An}{N}+P(\frac{n}{N})\bigr)
		\lesssim N (qD\tau)^{-1/2} \log(2N).
	\end{equation}
\end{enumerate}
\end{lemma}

\begin{proof}
For (1), let $I_w\subset I_{\mathrm{reg}}$ have length $O(\Delta)$ and contain $\operatorname{supp}w$. Applying \eqref{eq:poisson} we have the left-hand side of \eqref{eq:smooth-cap}
$$
=  \frac{N}{q} \sum_{m \in \mathbb{Z}} G_q(m,c) I_m(w) = S_1 + S_2, 
$$
where
$$
S_1 = \frac{N}{q} \sum_{\substack{m \in \mathbb{Z} \\ \operatorname{dist}(m\lambda-A,P'(I_w))<\lambda}} G_q(m,c) I_m(w)
$$
and
$$
S_2 = \frac{N}{q} \sum_{\substack{m \in \mathbb{Z} \\ \operatorname{dist}(m\lambda-A,P'(I_w)) \ge \lambda}} G_q(m,c) I_m(w)
$$
For $S_1$, we have $|I_w|\sup_{I_w}|P''|\lesssim\Delta D\tau=\lambda$, so the number of such $m$ is $O(1)$. By the second-derivative test it follows that 
$$
S_1 \lesssim N q^{-1/2} B_w(D\tau)^{-1/2}. 
$$
Applying Lemma~\ref{lem:nonstationary} we derive that
$$
S_2 \lesssim \frac{q^{1/2} B_w}{\lambda \Delta} \lesssim N (qD\tau)^{-1/2} B_w.
$$

For (2), take a smooth bump $\chi_I$ equal to $1$ on $I$, supported in the $\Delta$-neighborhood $I'$ of $I$, with
$$
\max_{0\le r\le2}\Delta^r\|\chi_I^{(r)}\|_\infty\lesssim1.
$$
 Then  $w=\eta\chi_I\in C_c^\infty(I_{\mathrm{reg}})$ satisfies $B_w\lesssim1$ (since $\Delta\lesssim\tau$), by \eqref{eq:smooth-cap} the smoothed sum
 $$
 \sum_{n\in\mathbb Z} w\big( \frac{n}{N} \big) \mathrm{e}\bigl(\frac{cn^2}{q}+\frac{An}{N}+P(\frac{n}{N})\bigr) \lesssim N (qD\tau)^{-1/2}.
 $$
uniformly in $A$. We finish the proof by appling \eqref{eq:completion} to the left-hand side of \eqref{eq:sharp-cap} with $p=\infty$, 
$$
a_n = \mathbf{1}_{I}(\frac{n}{N}), \qquad b_n = w (\frac{n}{N}) \mathrm{e} \big(\frac{cn^2}{q} + P(\frac{n}{N}) \big).
$$ 
Here we expand $a_n$ on $\mathbb Z/P_*\mathbb Z$ to get $C\lesssim\log(2N)$, where $P_*\asymp N$ is larger than the support length of $w$. 
\end{proof}

Before the next lemma, we define some terms under the hypotheses of Lemma~\ref{lem:decomposition} with $H > H_0$. Fix a small constant $c_1>0$ so that 
$$
K_0=[3/4-2c_1,3/4+2c_1]\subset(1/2,1), 
$$
with
$$
3c_1\tau<\min(1/4-t_-,t_+-5/4) \quad \text{~and~} \quad 4c_1\tau<1/3
$$
for all $\tau\lesssim1$. Fix a smooth cutoff $0\le\gamma_0\le1$ in $C_c^\infty(K_0)$ equal to $1$ on $[3/4-c_1,3/4+c_1]$. Recall that $0<V\le U$. Define that
$$
U_D=\min(U,D), \quad V_D=\min(V,D)
$$
and
$$
L=D\tau^3, \quad L_3=\tau^3 U_D, \quad  W=\tau^4 V_D.
$$
Let 
$$
\mathcal A\subset\{(B,u,v)\in\mathbb R\times[0,U]\times[0,V]: 1+|B|+u+v\asymp D\}
$$
be measurable. 

Choose a smooth partition $\{\chi_{t_0}\}$ on $(1/4, 5/4)$ with $t_0$ spaced by $c_1 \tau/2$ such that
$$
\sum_{t_0} \chi_{t_0} (t) = 1, \quad \|\partial_t^r\chi_{t_0}\|_\infty\lesssim_r\tau^{-r}
$$
and
$$
\operatorname{supp} \chi_{t_0} \subset [t_0-c_1\tau,t_0+c_1\tau],
$$
which meets $(1/4,5/4)$. Note that the number of such $t_0$ is $O(\tau^{-1})$. For each $t_0$, set 
$$
h=t_0-\frac{3}{4}\tau
$$
and
\[
\mathcal A_\tau(t_0)=\{(B,u,v)\in\mathcal A:
\operatorname{supp}(\chi_{t_0}\psi_\tau(B,u,v,\cdot))\ne\emptyset\}.
\]
Define $\Psi_{t_0}:(A,B,u,v)\mapsto X=(X_1,X_2,X_3,X_4)$ by
\[
X_j=\tau^j\sum_{\ell=j}^4\binom{\ell}{j}h^{\ell-j}x_\ell
\qquad \bigl((x_1,x_2,x_3,x_4)=(A,B,u,v)\bigr).
\]
Write $\Phi_{t_0}:(B,u,v)\mapsto(X_2,X_3,X_4)$ and put
\begin{equation}\label{eq:F0}
F_{0,t_0}(X)=\sum_{n\in\mathbb Z}\gamma_0(r_n)\mathrm{e}(cn^2q^{-1})
\mathrm{e}(X\cdot\Gamma(r_n))
\end{equation}
with 
$$
r_n=\frac{n/N-h}{\tau}. 
$$

Now we are ready for the next lemma. 
\begin{lemma}\label{lem:band-reduction}
We have the following results.
\begin{enumerate}
	\item $\det(\partial X/\partial(A,B,u,v))=\tau^{10}.$
	\item $\displaystyle x\cdot\Gamma(h+\tau r)=\sum_{\ell=1}^4 x_\ell h^\ell+X\cdot\Gamma(r).$
	\item $h+\tau K_0\subset[t_-,t_+]$ with $|h+\tau K_0|=4c_1\tau<1/3$. 
	\item $P''(h+\tau r)=\tau^{-2}(2X_2+6X_3r+12X_4r^2).$
	\item $P'''(h+\tau r)=\tau^{-3}(6X_3+24X_4r).$
	\item For sufficiently small fixed $c_1>0$, every $(B,u,v)\in\mathcal A_\tau(t_0)$ satisfies $|P''(h+\tau r)|\asymp D\tau$ and $|P'''(h+\tau r)|\lesssim D$ for all $r\in K_0$, 
	\item There is an absolute constant $C>0$ such that
	\begin{equation}\label{eq:box-enclosure}
		\Phi_{t_0}(\mathcal A_\tau(t_0))
		\subset[-CL,CL]\times[-CL_3,CL_3]\times[0,CW].
	\end{equation}

\end{enumerate}
\end{lemma}
\begin{proof}
(1) is proved by the triangular structure of $\Psi_{t_0}$. 

(2)(4)(5) Expanding $x\cdot\Gamma(h+\tau r)$ in powers of $r$ the proof is done.  

(3) Since $[t_0-c_1\tau,t_0+c_1\tau]$ meets
$(1/4,5/4)$, we have $h+\tau K_0=[t_0-2c_1\tau,t_0+2c_1\tau]\subset[t_-,t_+]$
and $|h|\lesssim1$. 

(6) If $(B,u,v)\in\mathcal A_\tau(t_0)$, pick
$t_1\in\operatorname{supp}(\chi_{t_0}\psi_\tau)\subset[t_0-c_1\tau,t_0+c_1\tau]$;
then $|P''(t_1)|\asymp D\tau$ and $|P'''|\lesssim D$ on $[t_-,t_+]$, so
taking $c_1$ small preserves $|P''(h+\tau r)|\asymp D\tau$ for all
$r\in K_0$. 

(7) $0\le u\lesssim U_D$ and $0\le v\lesssim V_D\le U_D$
give $0\le X_4=\tau^4v\lesssim W$ and
$$
|X_3|=\tau^3|u+4hv|\lesssim\tau^3U_D=L_3\le L.
$$
Apply (4) with $r=3/4$ yields
$$
2X_2=\tau^2P''(t_0)-\frac{9}{2}X_3-9X_4,
$$
whence $|X_2|\lesssim L$. 
\end{proof}

\begin{lemma}\label{lem:2.5.2}
	For every $1\le p<\infty$,
	\begin{equation}\label{eq:band-reduction}
		\int_{\mathcal A}\int_0^N|T_\tau(A)|^p\,dA\,dB\,du\,dv
		\lesssim_p \frac{N\tau^{-9}}{L}\sum_{t_0}\sup_{s\in\mathbb R}
		\int_{[s,s+L]\times\Phi_{t_0}(\mathcal A_\tau(t_0))}
		|F_{0,t_0}(X)|^p\,dX.
	\end{equation}
\end{lemma}
\begin{proof}
Define
$$
T_{\tau,t_0}(A)=\sum_{n\in\mathbb Z}\chi_{t_0}(\frac{n}{N}) \psi_\tau(B,u,v; \frac{n}{N}) \mathrm{e}\big(\frac{cn^2}{q}+A\frac{n}{N}+P(\frac{n}{N})\big).
$$
By Lemma~\ref{lem:decomposition}, for each $(B,u,v)\in\mathcal A$, $T_{\tau,t_0}(A) \neq 0$ for at most $O(1)$ times, so
$$
|T_\tau(A)|^p\lesssim_p\sum_{t_0}\mathbf{1}_{\mathcal A_\tau(t_0)} |T_{\tau,t_0}(A)|^p. 
$$
On $\mathcal A_\tau(t_0)$,
$$
\eta_{t_0}(r)=\chi_{t_0}(h+\tau r)\psi_\tau(B,u,v;h+\tau r)
$$
is supported in $[3/4-c_1,3/4+c_1]\subset(1/2,1)$, satisfying
$\chi_{t_0}\psi_\tau=\gamma_0\eta_{t_0}$, and has
$$
\|\eta_{t_0}\|_\infty+\|\eta_{t_0}''\|_\infty\lesssim1.
$$
Expand $\eta_{t_0}$ in a Fourier series with 
$$
\eta_{t_0} = \sum_{k \in \mathbb{Z}} d_k \mathrm{e}(kr),
$$
where
$$
d_k \lesssim (1+ |k|)^{-2}, \qquad \sum_k |d_k| \lesssim 1. 
$$
Write 
$$
G_{t_0} (A; B, u, v) = \sum_{n\in\mathbb Z} \gamma_0(r_n)  \mathrm{e}\big(\frac{cn^2}{q}+A\frac{n}{N}+P(\frac{n}{N})\big).
$$
By \eqref{eq:completion} we have that
$$
\|T_{\tau, t_0}(A)\|_{L^p(0,N)}\ \lesssim \|G_{t_0}(\cdot; B, u, v)\|_{L^p(0,N)}. 
$$
Changing variables via $\Psi_{t_0}$, since $\vert F_{0,t_0}(X_1+N\tau,X')\vert=\vert F_{0,t_0}(X_1,X')\vert$,
\begin{align*}
	\int_{\mathcal A}\int_0^N|T_\tau(A)|^p\,dA\,dB\,du\,dv&\lesssim_p\tau^{-10} \sum_{t_0} \int_{\Phi_{t_0}(\mathcal{A}_\tau(t_0))}\int_0^{N\tau}\vert F_{0,t_0}(X_1,X')\vert^pdX_1dX'\\
	&\le \tau^{-10} \sum_{t_0} \int_{\Phi_{t_0}(\mathcal{A}_\tau(t_0))}\frac{N\tau}{L}\sup_{s\in\mathbb{R}}\int_s^{s+L}\vert F_{0,t_0}(X_1,X')\vert^pdX_1dX'
\end{align*}
where $X' =(X_2,X_3,X_4)$.
\end{proof}

\section{Proof of Theorem~\ref{thm:main}}


For $R\ge2$, $J_*\ge1$, we define $M_4(J_*,R)$ as
$$
\begin{aligned}
\{\xi\in\mathbb R^4:\;&\xi_1\in[\tfrac12,1],\quad |\xi_2-\xi_1^2|\le J_*^{-2}, |\xi_3-3\xi_1\xi_2+2\xi_1^3|\le R^{-1},\\
&|\xi_4+4\xi_1\xi_3-9\xi_2^2+4\xi_1^4|\le R^{-1}\}.
\end{aligned}
$$
In application, we will take $J_*=R^\beta$ with $\frac13\le\beta\le\frac12$.

We will consider Schwartz functions $F$ such that
$\operatorname{spt}\widehat F\subset M_4(J_*,R)$.
Let $K_0\Subset(\frac12,1)$. We observe that
$\Gamma(K_0)+O(R^{-1})\subset M_4(J_*,R)$ since $R^{-1}\le J_*^{-2}$.

Let us divide $M_4(J_*,R)$ into boxes whose first coordinate are of length $J_*^{-1}$, i.e.
$$
M_4(J_*,R)\cap\{jJ_*^{-1}\le\xi_1<(j+1)J_*^{-1}\},\quad j\in\mathbb Z.
$$
We denote such a box by $\gamma$.

Let $\alpha\in(0,\beta]$, and $U=R^\alpha$. Fix $p\in(10,12]$, and set
$$
k=\max\{\tfrac52,\tfrac{p-6}2\}\in[\tfrac52,3].
$$
We assume further that
$$
RU\ge J_*^k.
$$
In the physical space, we consider boxes $\mathcal B$ of dimension $R\times R\times R\times U$.
It follows from \cite[Thm.~1.2]{Glidewell} and H\"older's inequality that, for any $\delta>0$,
\begin{equation}\label{eq:QuarticWeylSum:small-cap-input}
\int_{\mathcal B}|F|^p\lesssim_{\alpha,\beta,\delta}
R^\delta\left(\sup_\gamma\|F_\gamma\|_{L^\infty(\mathbb R^4)}\right)^p
J_*^{p/2}|\mathcal B|,
\end{equation}
where $\operatorname{spt}\widehat F\subset\Gamma(K_0)+O(R^{-1})$, and
$\widehat F_\gamma=\widehat F|_\gamma$.

\begin{lemma}\label{lem:QuarticWeylSum:localized-estimate}
Fix $10<p\le12$ and $k=\max\{\frac52,\frac{p-6}2\}$. Suppose
$$
0<L\lesssim N^2,\quad 1\le J\lesssim L^{1/2},\quad \Lambda\gtrsim N^{-1}.
$$
Let $F_0(X):=\sum_{\nu=1}^m d_\nu e(X\cdot\Gamma(s_\nu))$, where $m\le N$,
$|d_\nu|\le1$, and $s_\nu\in K_0$. Let $Q_-$ and $Q_+$ be
concentric cubes of side lengths $\sim L$ such that $Q_-\subset Q_+$ and 
$\operatorname{dist}(\partial Q_-,\partial Q_+)\sim L$.
Assume that for any interval $I$ of length at most $J^{-1}$, we have
\begin{equation}\label{eq:QuarticWeylSum:local-cap-hypothesis}
\sup_{X\in Q_+}\left|\sum_{s_\nu\in I}d_\nu e(X\cdot\Gamma(s_\nu))\right|\lesssim\Lambda.
\end{equation}
If $\mathcal Q\subset Q_-$ is a box of dimension $\sim L\times L\times L\times W$
with $0<W\lesssim L$, then for any $\varepsilon>0$,
\begin{equation}\label{eq:QuarticWeylSum:localized-estimate}
\int_{\mathcal Q}|F_0|^p\lesssim_\varepsilon
N^\varepsilon\Lambda^pJ^{p/2}L^3\max\{W,1,J^k/L\},
\end{equation}
where the implicit constant is independent of $N,L,J,\Lambda,\{s_\nu\}$, and $\{d_\nu\}$.

\end{lemma}

\begin{proof}
By \eqref{eq:QuarticWeylSum:local-cap-hypothesis}, it is trivial to have \eqref{eq:QuarticWeylSum:localized-estimate} when $N$ is finite or $J\le N^{\varepsilon/10}$.
We assume below that $N$ is large and $J\ge N^{\varepsilon/10}$.

Let $\ell=LN^{-\delta}$ with $\delta=\varepsilon/20$, and $R=C\min(\ell,J^3)$,
where $C=C(K_0)\ll1$. As $J\lesssim L^{1/2}$ and $J>N^{2\delta}$, we
obtain $L\gtrsim J^2>N^{4\delta}$. As a result,
$J^2N^{-\delta}\lesssim R\le CJ^3$, which implies further that
$2\lesssim N^{3\delta}\lesssim R\lesssim N^2$, and $\log J/\log R\ge\frac13$.

Divide $[\frac13,\frac12]=\bigcup I_j$ with $I_j=[j\delta,(j+1)\delta]$, and
take $\beta$ as the largest number of the form $j\delta$ such that
$$
\beta\le\min(\log J/\log R,\tfrac12).
$$
Set
$$
J_*=R^\beta,\quad \alpha=\max(\delta,k\beta-1),\quad U=R^\alpha.
$$
Then
$$
R^{1/3}\le J_*\le\min(J,R^{1/2}),\quad J/J_*\lesssim N^{2\delta}.
$$
Moreover, one can check that
$$
U\lesssim N^{2\delta}\max(1,J^k/L).
$$

Next we partition $\mathcal Q$ into boxes $E$ of dimension
$\sim\ell\times\ell\times\ell\times w_E$ with $w_E\le\ell$.
We denote the center of $E$ by $a_E$.
Let $\psi$ be a nonnegative function on $\mathbb R^4$ such that
$\operatorname{spt}\widehat\psi\subset B(0,1)$, and $|\psi|\ge1$ on $[-1,1]^4$.
Denote
$$
F_E(X)=\psi((X-a_E)/\ell)F_0(X).
$$
Then $\operatorname{spt}\widehat F_E\subset\Gamma(K_0)+O(\ell^{-1})$.
We recall that $\ell^{-1}\lesssim R^{-1}$. Therefore we can apply \eqref{eq:QuarticWeylSum:small-cap-input}.
We will estimate $\|(F_E)_\gamma\|_{L^\infty}$ below.

Firstly, we divide $K_0$ into subintervals $\kappa$ of length $J_*^{-1}$, and define
$$
F_{0,\kappa}(X)=\sum_{s_\nu\in\kappa}d_\nu e(X\cdot\Gamma(s_\nu)).
$$
It follows from \eqref{eq:QuarticWeylSum:local-cap-hypothesis} and $|\kappa|=J_*^{-1}$ that
\begin{equation}\label{eq:QuarticWeylSum:localized-cap-decay}
\begin{aligned}
|\psi((X-a_E)/\ell)F_{0,\kappa}(X)|
&\lesssim_\delta\left(\frac J{J_*}\right)\Lambda(1+|X-a_E|/\ell)^{-2}
\end{aligned}
\end{equation}
when $X\in Q_+$. For $X\in(Q_+)^c$, using the rapid decay
of $\psi$, $|X-a_E|/\ell\gtrsim N^\delta$, and the trivial estimate
$$
|F_{0,\kappa}(X)|\lesssim m\lesssim N^2\Lambda,
$$
we establish \eqref{eq:QuarticWeylSum:localized-cap-decay} as well.

Let us write $y=X_1-(a_E)_1$.
Let $I_\gamma$ be the first-coordinate interval of  $\gamma$
Applying \eqref{eq:QuarticWeylSum:localized-cap-decay}, we see that
\begin{equation}\label{eq:QuarticWeylSum:global-cap-bound}
\begin{aligned}
|(F_E)_\gamma(X)|
&\lesssim \int|(\mathbf 1_{I_\gamma})^\vee(t)|\frac J{J_*}\Lambda
\left(1+\frac{|(y-t,(X-a_E)')|}{\ell}\right)^{-2}\,dt\\
&\lesssim\frac J{J_*}\Lambda\int\min(J_*^{-1},|t|^{-1})
\left(1+\frac{|y-t|}{\ell}\right)^{-2}\,dt\\
&\lesssim\frac J{J_*}\Lambda\left(\int_{|t|\le J_*}J_*^{-1}\,dt
+\int_{J_*<|t|\le\ell}|t|^{-1}\,dt+\ell^{-1}\int_{|t|>\ell}
\left(1+\frac{|y-t|}{\ell}\right)^{-2}\,dt\right)\\
&\lesssim\frac J{J_*}\Lambda\log(2N),
\end{aligned}
\end{equation}
where we use
$$
|(\mathbf 1_{I_\gamma})^\vee(t)|\lesssim\min(J_*^{-1},|t|^{-1}),
$$
and the notation $X'=(X_2,X_3,X_4)$.

As each $E$ is covered by $O((\frac\ell R)^3(1+\frac{w_E}U))$
boxes $\mathcal B$ of dimension $R\times R\times R\times U$, we apply \eqref{eq:QuarticWeylSum:small-cap-input}
and \eqref{eq:QuarticWeylSum:global-cap-bound} to obtain
$$
\begin{aligned}
\int_E|F_0|^p&\lesssim\int_E|F_E|^p \lesssim_\delta R^\delta\Lambda^p\left(\frac J{J_*}\right)^p
J_*^{p/2}[\log(2N)]^p\ell^3\max(w_E,U).
\end{aligned}
$$
Covering $\mathcal Q$ by $E$, we obtain the control
$$
\begin{aligned}
\int_{\mathcal Q}|F_0|^p
&\lesssim_\delta\left(\frac L\ell\right)^3\max(W,U)R^\delta\Lambda^p
\left(\frac J{J_*}\right)^p(J_*)^{p/2}[\log(2N)]^p\ell^3\\
&\lesssim_\delta L^3N^\varepsilon\max(W,1,J^k/L)\Lambda^pJ^{p/2}.
\end{aligned}
$$
\end{proof}

\begin{proposition}\label{prop:QuarticWeylSum:common-range}
Let $10<p\le12$, $k=\max\{\frac52,\frac{p-6}2\}$, and suppose
$a+b=\frac{p-6}2$, $a\ge b\ge0$, and $b\ge a(k-2)/2$.
Then
$$
\int_{\Omega_N(a,b)}|S_{I,N}(x)|^p\,dx\lesssim_{p,a,\varepsilon}N^{p+\varepsilon}
$$
for all $ N\ge2$, and $I\subset[\frac N2,N]$.

\end{proposition}

\begin{proof}
Recalling that
\begin{equation}\label{eq:QuarticWeylSum:interval-projection}
\int_{\Omega_N(a,b)}|S_{I,N}(x)|^p\,dx
\lesssim\int_{\Omega_N(a,b)}|S_\gamma|^p\,dx,
\end{equation}
it suffices to control the RHS.

Set $U=N^a$, $V=N^b$, and $Q=\lfloor N^{1-\frac a2}\rfloor$. For $q\le Q$, we define
$R_q=\frac{N^2}{qQ}$, and
$$
\operatorname{arc}(\tfrac cq):=[0,N]\times[-R_q,R_q]\times[0,U]\times[0,V],
$$
where $1\le V\le U\le\frac{N^2}{Q^2}$.
It follows from Dirichlet's approximation theorem and the $N^2$-periodicity of $S_\gamma$
in $x_2$ that
\begin{equation}\label{eq:scales}
\int_{\Omega_N(a,b)}|S_\gamma(x)|^p\,dx \le\sum_{q=1}^Q\sum_{\substack{c\bmod q\\(c,q)=1}}
\int_{\operatorname{arc}(c/q)}
|S_\gamma(A,\tfrac{N^2c}q+B,u,v)|^p\,dA\,dB\,du\,dv.
\end{equation}
Applying \eqref{eq:exceptional}, it is further bounded by
$$
N^p+\sum_{q=1}^Q\sum_{\substack{c\bmod q\\(c,q)=1}}
\int_{\operatorname{arc}(c/q)}
|\mathcal R(A,\tfrac{N^2c}q+B,u,v)|^p\,dA\,dB\,du\,dv.
$$

Let us fix $q\le Q$ and $c$ with $(c,q)=1$. On $[-R_q,R_q]\times[0,U]\times[0,V]$, we have
$$1+|B|+u+v\le4R_q\le4\frac{N^2}{q^2}.$$
Partition $[-R_q,R_q]\times[0,U]\times[0,V]$ into dyadic pieces $\mathcal D_q(D)$
with $1+|B|+u+v\sim D$ in $\mathcal D_q(D)$, where $D\le4N^2/q^2$.
Set $U_D=\min(N^a,D)$, and $V_D=\min(N^b,D)$, then $1\le V_D\le U_D\le D$,
$|\mathcal D_q(D)|\lesssim DU_DV_D\lesssim D^2V_D$, and
\begin{equation}\label{eq:QuarticWeylSum:dyadic-sum}
\sum_{D\lesssim R_q}D^2V_D\lesssim R_q^2N^b.
\end{equation}

We are able to apply Lemma~\ref{lem:decomposition}.
If $H=D\frac qN\le H_0$, we have
\begin{equation}\label{eq:QuarticWeylSum:small-H-remainder}
|\mathcal R|\lesssim\sqrt N.
\end{equation}
If $H>H_0$, we have
\begin{equation}\label{eq:QuarticWeylSum:band-remainder}
\left|\mathcal R-\sum_\tau T_\tau\right|\lesssim\sqrt N.
\end{equation}
We will estimate $T_\tau$ below when $H>H_0$.

We will apply Lemma~\ref{lem:band-reduction} with $\mathcal A=\mathcal D_q(D)$.
Fix $t_0$ and $s$, and put
$\mathcal B_{t_0}(s)=[s,s+L]\times\Phi_{t_0}(\mathcal A_\tau(t_0))$.
It follows from Lemma~\ref{lem:band-reduction} that
$$
\mathcal B_{t_0}(s)\subset[s,s+L]\times[-CL,CL]^2\times[0,CW].
$$
Let $\delta_0>0$ be small enough, and cover $\mathcal B_{t_0}(s)$ by boxes
$\mathcal Q$ of dimension $\sim\delta_0L\times\delta_0L\times\delta_0L\times w_{\mathcal Q}$
with $w_{\mathcal Q}\le\min(CW,\delta_0L)$.

Fix $X^*\in\mathcal Q\cap\mathcal B_{t_0}(s)$.
Set $Q_-:=X^*+[-\delta_0L,\delta_0L]^4$, and
$Q_+:=X^*+[-2\delta_0L,2\delta_0L]^4$.
Since $(X_2^*,X_3^*,X_4^*)\in\Phi_{t_0}(\mathcal A_\tau(t_0))$,
by Lemma~\ref{lem:band-reduction}, the polynomial $P_*$ related to $\Psi_{t_0}^{-1}(X^*)$ satisfies
$$
\begin{aligned}
|P_*''(h+\tau r)|&\sim D\tau,\\
|P_*'''(h+\tau r)|&\lesssim D
\end{aligned}
$$
on $K_0$. Moreover, using that $|X_j-X_j^*|\le2\delta_0L$, we have
$
|P''|\sim D\tau
$ and 
$|P'''|\lesssim D
$
throughout $h+\tau K_0$ for any $ X\in Q_+$.

We apply Lemma~\ref{lem:QuarticWeylSum:localized-estimate} below. Set
$$
L=D\tau^3,\quad W=\tau^4V_D,\quad J=\frac{D\tau^2q}N,\quad
M=\frac N{\sqrt{qD\tau}}.
$$
One can check that $J^k\lesssim LV_D$.

By Lemma~\ref{lem:cap}, \eqref{eq:QuarticWeylSum:local-cap-hypothesis} is satisfied with $\Lambda=M\log(2N)$.
As a result, by Lemma~\ref{lem:QuarticWeylSum:localized-estimate},
$$
\int_{\mathcal B_{t_0}(s)}|F_{0,t_0}(X)|^p\,dX
\lesssim_\varepsilon N^\varepsilon M^pJ^{p/2}L^3V_D.
$$

Recalling Lemma~\ref{lem:2.5.2}, we obtain further that
$$
\begin{aligned}
&\int_{\mathcal D_q(D)}\int_0^N|T_\tau|^p\,dA\,dB\,du\,dv\lesssim_\varepsilon N^{\frac p2+1+\varepsilon}D^2V_D\tau^{\frac p2-4}.
\end{aligned}
$$
As $\frac p2-4>1$, by Minkowski's inequality, we obtain from \eqref{eq:QuarticWeylSum:small-H-remainder} and \eqref{eq:QuarticWeylSum:band-remainder} that
$$
\int_{\mathcal D_q(D)}\int_0^N|\mathcal R|^p
\lesssim_\varepsilon N^{\frac p2+1+\varepsilon}D^2V_D.
$$
Summing over $D\lesssim R_q$, $c$ with $(c,q)=1$, and $1\le q\le Q$, we obtain
$$
\begin{aligned}
\sum_{q=1}^Q\sum_{\substack{c\bmod q\\(c,q)=1}}
\int_{\operatorname{arc}(c/q)}|\mathcal R|^p
&\lesssim_\varepsilon N^{\frac p2+1+\varepsilon}\sum_{q\le Q}qR_q^2V
\lesssim_\varepsilon N^{p+\varepsilon},
\end{aligned}
$$
where we use \eqref{eq:QuarticWeylSum:dyadic-sum}, $Q\sim N^{1-\frac a2}$, and
$$
\sum_{q=1}^Qq^{-1}\lesssim\log(2N).
$$
This completes the proof by \eqref{eq:QuarticWeylSum:interval-projection}.\end{proof}


We are ready to prove Theorem~\ref{thm:main}.
\begin{proof}[Proof of Theorem~\ref{thm:main}]
For $a\le2$, we have $k=a+b$, which implies that
$b-\frac{a(k-2)}2=k(1-\frac a2)\ge0$. Thus we can apply Proposition~\ref{prop:QuarticWeylSum:common-range} in this case.

For $a\ge2$, and $b>0$, we have $a+b\in[\frac52,3]$.
As a result, $2\le a\le3$, and $0<b\le1$. Then we apply \cite[Cor.~1.4]{Glidewell} to obtain the estimate.

For the endpoint case $b=0$, we have $a=(p-6)/2\in[5/2,3]$.
Take $\delta>0$ small enough, and define $p'=p+2\delta<14$.
Let $a'=a$, $b'=b+\delta$.
Then $a'+b'=\frac{p'-6}2$, and $\Omega_N(a,b)\subset\Omega_N(a',b')$.
Applying \cite[Cor.~1.4]{Glidewell}, we see that
$$
\begin{aligned}
\int_{\Omega_N(a,b)}|S_{I,N}(x)|^p\,dx
&\lesssim|\Omega_N(a,b)|^{1-\frac p{p'}}
\left(\int_{\Omega_N(a',b')}|S_{I,N}(x)|^{p'}\,dx\right)^{p/p'}\lesssim N^{p+2\delta}.
\end{aligned}
$$
\end{proof}

\section{Proof of Theorem~\ref{thm:tenth}}

Put $\Gamma_3(r)=(r,r^2,r^3)$. For an integer $j \ge 2$, define 
\[
\mathcal M_j=\bigl\{\xi\in\mathbb R^3:\ 0\le\xi_1\le1,\quad
|\xi_2-\xi_1^2|\le j^{-2},\quad
|\xi_3-3\xi_1\xi_2+2\xi_1^3|\le j^{-2}\bigr\}.
\]
Applying \cite[Theorem 3]{GuthMaldague} with 
$$
R=j^2, \beta=1/2, p=10
$$
gives, for every $\eta>0$ and every
Schwartz function $f$ on $\mathbb R^3$ with
$\operatorname{supp}\widehat f\subset\mathcal M_j$,
\begin{align}\label{eq:cubic-input}
	\nonumber \|f\|_{L^{10}(\mathbb R^3)}^{10} &\lesssim_\eta R^{\eta} (R^{\beta(p/2-1)} + R^{\beta(p-4)-1})  \sum_\gamma\|f_\gamma\|_{L^{10}(\mathbb R^3)}^{10} \\
	&\lesssim_\eta j^{4+\eta}\sum_\gamma\|f_\gamma\|_{L^{10}(\mathbb R^3)}^{10},
\end{align}
where $f_\gamma$ is the sharp Fourier restriction to the cap 
$$
\gamma = \{ \xi \in \mathcal M_j: lj^{-1} \le \xi_1 < (l+1)j^{-1} \}. 
$$
\begin{lemma}\label{lem:tenth-cube}
	Under the hypotheses of Lemma~\ref{lem:cap}, let $F_{0,t_0}(\mathbf x,X_4)$ be defined as \eqref{eq:F0} where
	$\mathbf x=(X_1,X_2,X_3)$. Fix $X_4\in\mathbb R$ and a center
	$\mathbf x_c\in\mathbb R^3$ such that the polynomial $P_c$ corresponding
	to $(\mathbf x_c,X_4)$ under $\Psi_{t_0}^{-1}$ satisfies
	\[
	|P_c''(h+\tau r)|\asymp D\tau,\qquad |P_c'''(h+\tau r)|\lesssim D
	\qquad (r\in K_0).
	\]
	Set $j=\max(2,\lceil J\rceil)$ and $\rho=j^2$. Then for every $\eta>0$,
	on the cube $Q_\rho(\mathbf x_c)=\mathbf x_c+[-\rho,\rho]^3$,
	\begin{equation}\label{eq:tenth-cube}
		\int_{Q_\rho(\mathbf x_c)}|F_{0,t_0}(\mathbf x,X_4)|^{10}\,d\mathbf x
		\lesssim_\eta j^{5+\eta}M^{10}\rho^3.
	\end{equation}
\end{lemma}
\begin{proof}
Since $J\gtrsim1$, we have $j\asymp J$ and $\tau/j\asymp\Delta$. Choose
a smooth partition of unity $\{\chi_\kappa\}$ on $K_0$ with centers
$r_\kappa=\kappa/j$ satisfying
$$
\operatorname{supp}\chi_\kappa\subset[r_\kappa-j^{-1},r_\kappa+j^{-1}],
$$
and 
$$
\max_{0\le m\le2}j^{-m}\|\chi_\kappa^{(m)}\|_\infty\lesssim1.
$$
Note that $\kappa\in\mathbb Z$ and the number of indices whose supports meet $K_0$ is $O(j)$. 
Write $F_{0,t_0}=\sum_\kappa F_\kappa$, where
\[
F_\kappa(\mathbf x,X_4)
=\sum_{n\in\mathbb Z}\chi_\kappa(r_n)\gamma_0(r_n)\mathrm{e}(\frac{cn^2}{q}) \mathrm{e}\bigl(\mathbf x\cdot\Gamma_3(r_n)+X_4r_n^4\bigr).
\]
Put $\mathbf v=\mathbf x-\mathbf x_c=(v_1,v_2,v_3)$ and $z=j(r-r_\kappa)$ with $r = (t-h)/\tau$. Expanding $\mathbf v\cdot\Gamma_3(r)$ around $r_\kappa$ gives
\[
\mathbf v\cdot\Gamma_3(r)
=\mathbf v\cdot\Gamma_3(r_\kappa)
+\bigl(\mathbf v\cdot\Gamma_3'(r_\kappa)\bigr)(r-r_\kappa)
+\Psi_\kappa(z),
\]
where
\[
\Psi_\kappa(z)=\frac{v_2+3r_\kappa v_3}{j^2}\,z^2+\frac{v_3}{j^3}\,z^3.
\]
Take 
$$
A' = A_c + \frac{\mathbf{v} \cdot \Gamma_3'(r_\kappa)}{\tau}. 
$$
Hence
\[
|F_\kappa(\mathbf x,X_4)|
=\left|\sum_{n\in\mathbb Z}w_\kappa(n/N)
\mathrm{e}\bigl(cn^2/q+A'(n/N)+P_c(n/N)\bigr)\right|,
\]
where
$$
w_\kappa(t)=\chi_\kappa(r)\gamma_0(r)\mathrm{e}(\Psi_\kappa(z)). 
$$
Since $\operatorname{supp}w_\kappa$ lies in an
interval of length $O(\tau/j)=O(\Delta)$ inside $h+\tau K_0$ where
$|z|\le1$, and $|\Psi_\kappa'(z)|,|\Psi_\kappa''(z)|\lesssim|\mathbf v|/\rho$
for $|z|\le1$, it gives
$$
B_{w_\kappa}=\max_{0\le m\le2}\Delta^m\|w_\kappa^{(m)}\|_\infty
\lesssim(1+|\mathbf v|/\rho)^2.
$$
Applying Lemma~\ref{lem:cap}(1) to $P_c$
yields
\[
|F_\kappa(\mathbf x,X_4)|
\lesssim M\Bigl(1+\frac{|\mathbf x-\mathbf x_c|}{\rho}\Bigr)^2
\qquad(\mathbf x\in\mathbb R^3).
\]
Fix a Schwartz function $\psi$ on $\mathbb R^3$ with
$\operatorname{supp}\widehat\psi\subset B(0,c)$ for small $c>0$
depending on $K_0$ and $|\psi|\ge1$ on $[-1,1]^3$. We put
$$
\psi_\rho(\mathbf x)=\psi((\mathbf x-\mathbf x_c)/\rho).
$$ 
Rapid decay of $\psi$ gives
$$
\|\psi_\rho F_\kappa(\cdot,X_4)\|_{L^{10}(\mathbb R^3)}^{10}
\lesssim M^{10}\rho^3.
$$ 
The Fourier support of $\psi_\rho F_{0,t_0}(\cdot,X_4)$
lies within $c/\rho=cj^{-2}$ of $\Gamma_3(K_0)$. Hence the support is inside
$\mathcal M_j$ for small $c$. Because the centers $r_\kappa=\kappa/j$
are spaced by $j^{-1}$ and $c/\rho\le j^{-1}$, each cap slab $\gamma$
of width $j^{-1}$ in $\xi_1$ meets $\operatorname{supp}\widehat{\psi_\rho F_\kappa}$
for only $O(1)$ indices $\kappa$. By boundedness of the one-dimensional
Riesz projection in $X_1$ on $L^{10}(\mathbb R^3)$
\cite{HollenbeckVerbitsky},
\[
\|(\psi_\rho F_{0,t_0})_\gamma\|_{L^{10}(\mathbb R^3)}^{10}
\lesssim\sum_{\kappa:\,\operatorname{supp}\widehat{\psi_\rho F_\kappa}\cap\gamma\ne\emptyset}
\|\psi_\rho F_\kappa\|_{L^{10}(\mathbb R^3)}^{10}
\lesssim M^{10}\rho^3.
\]
Since $|\psi_\rho|\ge1$ on $Q_\rho(\mathbf x_c)$ and $\mathcal M_j$ has
$j$ caps, applying \eqref{eq:cubic-input} to $\psi_\rho F_{0,t_0}(\cdot,X_4)$
proves \eqref{eq:tenth-cube}.
\end{proof}

Now we are ready for Theorem~\ref{thm:tenth}. 
\begin{proof}[Proof of Theorem~\ref{thm:tenth}]
Fix $1\le a\le2$ and $b=2-a$, and retain the rational-arc scales \eqref{eq:scales}. Similarly we have that
\[
\int_{\Omega_N(a,b)}|S_\gamma|^{10}\,dx
\le\sum_{q\le Q}\sum_{\substack{c\bmod q\\(c,q)=1}}
\int_{\operatorname{arc}(c/q)}|S_\gamma|^{10}.
\]
By \eqref{eq:tenth-near}, the exceptional part $S_{\mathrm{near}}$ contributes $\lesssim_\varepsilon N^{10+\varepsilon}$
uniformly in $a$.

It remains to bound $\mathcal R=S_\gamma-S_{\mathrm{near}}$ on each arc. By the same construction in the proof of Proposition \ref{prop:QuarticWeylSum:common-range}, the $10$th-power integral over $[0,N]\times\mathcal D_q(D)$ of the error term $O(N^{1/2})$ is $O(N^6DU_DV_D)$. Fix $D$ with $H>H_0$ and a signed dyadic band $T_\tau$ ($\tau_0\le\tau\lesssim1$). It is easy to check that $\rho \lesssim L_3 \le L$. 

Apply Lemma~\ref{lem:band-reduction} to $\mathcal A=\mathcal D_q(D)$. Fix $t_0$ and $s\in\mathbb R$, and put
$$
\mathcal B_{t_0}(s)=[s,s+L]\times\Phi_{t_0}(\mathcal A_\tau(t_0)). 
$$
By \eqref{eq:box-enclosure},
$$
\mathcal B_{t_0}(s)\subset[s,s+L]\times[-CL,CL]\times[-CL_3,CL_3]\times[0,CW].
$$
Freeze $X_4\in[0,CW]$ and partition 
$$
[s,s+L]\times[-CL,CL]\times[-CL_3,CL_3]
$$
into $O(L^2L_3/\rho^3)$ grid cubes of side $\rho$ (since $\rho\lesssim L_3\le L$). For each grid cube meeting the slice
$$
\mathcal S(X_4)=\{\mathbf x:(\mathbf x,X_4)\in\mathcal B_{t_0}(s)\},
$$
pick $\mathbf x_c\in\mathcal S(X_4)$ inside it. Then the grid cube lies
in $Q_\rho(\mathbf x_c)$ and the polynomial $P_c$ associated to
$(\mathbf x_c,X_4)$ satisfies $|P_c''|\asymp D\tau$ and
$|P_c'''|\lesssim D$ on $h+\tau K_0$. Applying Lemma~\ref{lem:tenth-cube}
with $\eta=\varepsilon$ on each nonempty cube, summing over the
$O(L^2L_3/\rho^3)$ cubes, and integrating over $X_4\in[0,CW]$ yields
(using $j^\varepsilon\lesssim_\varepsilon N^\varepsilon$)
\[
\int_{\mathcal B_{t_0}(s)}|F_{0,t_0}(X)|^{10}\,dX
\lesssim_\varepsilon N^\varepsilon M^{10}J^5L^2L_3W,
\]
uniformly in $a,s,t_0$. Substituting into \eqref{eq:band-reduction}
(with $O(\tau^{-1})$ centers $t_0$) gives
\begin{align}
	\int_{\mathcal D_q(D)}\int_0^N|T_\tau|^{10}\,dA\,dB\,du\,dv
	&\lesssim_\varepsilon N^\varepsilon\tau^{-11}\frac{N\tau}{L}
	M^{10}J^5L^2L_3W\notag\\
	&=N^{6+\varepsilon}DU_DV_D\tau^5.\label{eq:tenth-band}
\end{align}

Summing tenth roots of \eqref{eq:tenth-band} over dyadic $\tau\lesssim1$ and adding the $O(N^6DU_DV_D)$ error gives
$$
\int_{\mathcal D_q(D)}\int_0^N|\mathcal R|^{10}\lesssim_\varepsilon
N^{6+\varepsilon}DU_DV_D.
$$
Summing over $D\lesssim R_q$ gives
$$
\int_{\operatorname{arc}(c/q)}|\mathcal R|^{10}\lesssim_\varepsilon
N^{6+\varepsilon}R_qUV.
$$
Finally, summing over the $\varphi(q)$ reduced
numerators $c$ and $1\le q\le Q$ via
$$
\sum_{q\le Q}\varphi(q)R_q=(N^2/Q)\sum_{q\le Q}\varphi(q)/q\le N^2
$$
and $UV=N^2$ bounds the regular part by $O_\varepsilon(N^{10+\varepsilon})$.
Combining this with \eqref{eq:tenth-near} and applying
\eqref{eq:interval-projection} proves \eqref{eq:tenth-upper} uniformly
in $a\in[1,2]$.
\end{proof}	


\end{document}